\documentstyle[12pt]{article}
\def\lqq{\lq \lq }

\begin{document}
\rightline{November 1995}
\rightline{Preprint ICN-09-1995 ,  IFUNAM FT 95-82}

\begin{center}
{\LARGE Hidden $sl_2$-algebra of finite-difference equations}\footnote
{Talk presented at the Wigner Symposium, Guadalajara, Mexico, August 1995}
\vskip 0.5cm
Yuri SMIRNOV\footnote
{Instituto de F\'isica, UNAM, A.P. 20-364, 01000 M\'exico D.F.

On leave of absence from the Nuclear Physics Institute, Moscow State
University,
Moscow 119899, Russia\\
E-mail: smirnov@sysul2.ifisicacu.unam.mx}
and Alexander TURBINER\footnote
{Instituto de Ciencias Nucleares, UNAM, A.P. 70-543, 04510 M\'exico D.F.

On leave of absence from the Institute for Theoretical and Experimental
Physics,
Moscow 117259, Russia\\
E-mail: turbiner@nuclecu.unam.mx  or turbiner@axcrnb.cern.ch}
\end{center}
\vskip 0.3 cm
\begin{quote}
The connection between polynomial solutions of finite-difference equations
and finite-dimensional representations of the $sl_2$-algebra is established.
\end{quote}

\noindent
Recently it was found \cite{t1} that polynomial solutions of differential
equations
are connected to finite-dimensional representations of the
algebra $sl_2$ of first-order differential operators.
In this Talk it will be shown that there also exists a connection between
polynomial solutions of finite-difference equations
(like Hahn, Charlier and Meixner polynomials)
and unusual finite-dimensional representations of the algebra $sl_2$
of finite-difference operators. So, $sl_2$-algebra is the hidden algebra
of finite-difference equations with polynomial solutions.

First of all, we recall the fact that Heisenberg algebra
\begin{equation}
\label{e1}
[a,b] \equiv ab  -  ba \ =\ 1
\end{equation}
possesses several representations in terms of differential operators.
There is the standard, coordinate-momentum representation
\begin{equation}
\label{e2}
a\ =\ {d \over dx} \ ,\ b \ =\ x
\end{equation}
and recently  another one was found \cite{st}
\[
a\ =\ {\cal D}_+ \ ,
\]
\begin{equation}
\label{e3}
b\ =\ x(1 - \delta {\cal D}_-)\ ,
\end{equation}
where the ${\cal D}_{\pm}$ are translationally-covariant, finite-difference
operators
\begin{equation}
\label{e4}
{\cal D}_+ f(x) \ =\ {f(x+\delta) - f(x) \over \delta} \equiv
{(e^{\delta {d \over dx}} -1) \over \delta} f(x)
\end{equation}
and
\begin{equation}
\label{e5}
{\cal D}_- f(x) \ =\ {f(x) - f(x - \delta) \over \delta} \equiv
{ (1 - e^{ - \delta {d \over dx}}) \over \delta} f(x) \ ,
\end{equation}
where  $\delta$ is a parameter and ${\cal D}_+ \rightarrow {\cal D}_- $,
if $\delta \rightarrow -\delta$.

Now let us consider the Fock space over the operators $a$ and $b$ with a
vacuum $|0>$:
\begin{equation}
\label{e6}
a|0>\  = \ 0
\end{equation}
and define an operator spectral problem
\begin{equation}
\label{e7}
L[a, b] \varphi (b) = \lambda \varphi (b)
\end{equation}
where $L[\alpha, \beta]$ is a certain holomorphic function of the
variables $\alpha, \beta$. We will restrict ourselves studying the
operators $L[a, b]$ with polynomial eigenfunctions.

In \cite{t1} it was proven that $L$ has a certain number of
polynomial eigenfunctions if and only if, $L$ is the sum of two terms:
an element of the universal enveloping algebra of the $sl_2$-algebra taken
in the finite-dimensional irreducible representation
\begin{equation}
\label{e8}
J^+_n = b^2 a - n b,\ J^0_n = ba - {n \over 2},\ J^-_n=a
\end{equation}
where $n$ is a non-negative integer
\footnote{Taking $\ a,b \ $ from (2), the algebra (8) becomes the well-known
realization of $sl_2$ in first-order differential operators. If $a,b$ from (3)
are
chosen then (8) becomes a realization of $sl_2$ in finite-difference
operators.}
, and an annihilator $B(b)a^{n+1}$, where $B(b)$ is any operator function of
$b$.
The dimension of the representation (8) is equal
to $(n+1)$ and $(n+1)$ eigenfunctions of $L$ have the form of a polynomial
of degree not higher then $n$. These operators $L$ are named
{\it quasi-exactly-solvable}. Moreover,
if $L$ is presented as a finite-degree polynomial in the generators
$J^0 \equiv J^0_0$ and $ J^- \equiv J^-_0$ only, one
can prove that $L$ possesses infinitely-many polynomial eigenfunctions.
Such operators $L$ are named {\it exactly-solvable}.

It is evident that once the problem (7) is solved the eigenvalues will have
no dependence on the particular representation of  the operators $a$ and $b$.
This allows us to construct isospectral operators by simply taking different
representations of the operators $a$ and $b$ in the problem (7). In particular,
this implies that if we take the representation (3),  then the eigenvalues
of the problem (7) do not depend on the parameter $\delta$ !

Without loss of generality one can choose the vacuum
\begin{equation}
\label{e9}
|0>\  = \ 1
\end{equation}
and then it is easy to see that \cite{st}
\begin{equation}
\label{e10}
b^n |0> = [x(1 - \delta {\cal D}_-) ]^n |0> =\ x(x - \delta) (x - 2\delta)
\ldots (x - (n-1)\delta)
\equiv x^{(n)} \ .
\end{equation}
This relation leads to a very important conclusion: Once a solution of (7)
with $a,b$ (2) is found,
\begin{equation}
\label{e11}
\varphi (x) = \sum \alpha_k x^k\ ,
\end{equation}
then
\begin{equation}
\label{e12}
\tilde\varphi (x) = \sum \alpha_k x^{(k)}
\end{equation}
is the solution of (7) with $a,b$ given by (3).

Now let us proceed to a study of the second-order finite-difference equations
with
polynomial solutions and find the corresponding isospectral differential
equations.

The standard second-order finite-difference equation
relates an unknown function at three points and has the form\cite{nsu}
\begin{equation}
\label{e13}
 A(x) \varphi (x+\delta) - B(x)\varphi (x)+ C(x) \varphi (x-\delta)
  = \lambda \varphi (x),
\end{equation}
where $A(x), B(x), C(x)$ are arbitrary functions, $x \in R$.
One can pose a natural problem: {\it What are the most general coefficient
functions $A(x), B(x)$, $C(x)$ for which the equation (13) admits
infinitely-many
polynomial eigenfunctions ?} Basically, the answer is presented in \cite{t1}.
Any operator with the above property can be represented as a polynomial
in the generators $J^0, J^-$ of the $sl_2$-algebra:
\[
J^+= x({x \over \delta}-1) e^{ - \delta {d \over dx}} (1 - e^{ - \delta {d
\over dx}})\ ,
\]
\begin{equation}
\label{e14}
J^0= {x \over \delta} (1 - e^{ - \delta {d \over dx}})\ ,
\ J^-= {1 \over \delta} ( e^{  \delta {d \over dx}}-1)\ .
\end{equation}
which is the hidden algebra of our problem.
One can show that the most general polynomial in the generators (14)
leading to (13) is
\begin{equation}
\label{e15}
\tilde E = A_1 J^0J^0 (J^- + {1 \over \delta}) +A_2 J^0J^- +
 A_3 J^0 + A_4 J^- + A_5\ ,
\end{equation}
and in explicit form,
\[
[{A_4 \over \delta} + {A_2 \over \delta^2} x + {A_1 \over \delta^3} x^2]
e^{  \delta {d \over dx}} \ +
\]
\[
[A_5 - {A_4 \over \delta} + ({A_1 \over \delta^2} - 2 {A_2 \over \delta^2} \ +
{A_3 \over \delta})x - 2 {A_1 \over \delta^3} x^2] \ +
\]
\begin{equation}
\label{e16}
[ - ( {A_1 \over \delta^2} - {A_2 \over \delta^2} + {A_3 \over \delta}  )x\ +
{A_1 \over \delta^3} x^2] e^{  -\delta {d \over dx}}
\end{equation}
where the $A'$s are free parameters. The spectral problem corresponding to
the operator (16) is given by
\[
({A_4 \over \delta} + {A_2 \over \delta^2} x + {A_1 \over \delta^3} x^2)
f(x+\delta)\ -
\]
\[
[-A_5 + {A_4 \over \delta} - ({A_1 \over \delta^2} - 2 {A_2 \over \delta^2} +
{A_3 \over \delta})x + 2 {A_1 \over \delta^3} x^2] f(x) \ +
\]
\begin{equation}
\label{e17}
[ - ( {A_1 \over \delta^2} - {A_2 \over \delta^2} + {A_3 \over \delta}  )x+
{A_1 \over \delta^3} x^2]f(x-\delta) = \lambda f(x) \ .
\end{equation}
with the eigenvalues  $\lambda_k ={A_1 k^2 \over \delta} + A_3 k$.
This spectral problem has Hahn polynomials $h_k^{(\alpha,\beta)} (x, N)$
of the {\it discrete} argument $x=0,1,2\ldots, (N-1)$ as eigenfunctions (we use
the notation of \cite{nsu}). Namely,  these polynomials appear,
if $\delta = 1, A_5=0$ and
\[
A_1=-1,\ A_2= N-\beta-2,\ A_3= -\alpha-\beta-1,\ A_4=(\beta+1)(N-1) \ .
\]
If, however,
\[
A_1=1,\ A_2= 2-2N-\nu,\ A_3= 1-2N-\mu-\nu,\ A_4=(N+\nu-1)(N-1)
\]
the so-called analytically-continued Hanh polynomials $\tilde h_k^{(\mu,\nu)}
(x, N)$
 of the {\it discrete} argument $x=0,1,2\ldots, (N-1)$
appear, where $k=0,1,2\ldots$.

In general, our spectral problem (17) has Hahn polynomials
$h_k^{(\alpha,\beta)} (x, N)$
of the {\it continuous} argument $x$ as polynomial eigenfunctions. We must
emphasize
that
our Hahn polynomials of the {\it continuous} argument  {\it do not} coincide to
so-called
continuous Hahn polynomials known in literature\cite{hahn}.

So Equation (17) corresponds to the most general exactly-solvable
finite-difference problem, while the operator (15) is the most general element
of the universal enveloping $sl_2$-algebra leading to (13). Hence the Hahn
polynomials are related to the finite-dimensional representations of a certain
cubic element of the universal enveloping $sl_2$-algebra (for a general
discussion see \cite{t2}).

One can show that if the parameter $N$ is integer, then the higher Hahn
polynomials
$k\geq N$ have a representation
\begin{equation}
\label{e18}
h_k^{(\alpha,\beta)} (x, N)=x^{(N)}p_{k-N}(x)\  ,
\end{equation}
where $p_{k-N}(x)$ is a certain Hahn polynomial.
It explains an existence the only a finite number of the Hahn polynomials
of {\it discrete} argument $x=0,1,2\ldots, (N-1)$. Similar situation
occurs for the analytically-continued Hanh polynomials
\begin{equation}
\label{e19}
\tilde h_k^{(\mu,\nu)}  (x, N)=x^{(N)} \tilde p_{k-N}(x)\ ,\ k\geq N \ ,
\end{equation}
where $\tilde p_{k-N}(x)$ is a certain analytically-continued Hahn polynomial.

Furthermore, if we take the standard representation (8) for the algebra
$sl_2$ at $n=0$ with $a,b$ given by (2) and plug it into (15),
the third order differential operator {\it isospectral} to (16)
\begin{equation}
\label{e20}
\tilde E_2 ({d \over dx},x) =A_1 x^2 {d^3 \over dx^3} +
[(A_1+A_2)+ {A_1 \over \delta} x] x {d^2 \over dx^2} +
[A_4 + ({A_1 \over \delta} + A_3)x] {d \over dx} +A_5
\end{equation}
appears, which possesses polynomial eigenfunctions.

Taking in (17) $\delta=1, A_5=0$ and putting
\[
A_1=0, A_2= -\mu ,\ A_3= \mu -1, A_4= \gamma \mu \ ,
\]
and if $x=0,1,2\ldots, (N-1)$, we reproduce the equation, which has the
Meixner polynomials as eigenfunctions. Furthermore, if
\[
A_1=0, A_2= 0 ,\ A_3=  -1, A_4=  \mu \ ,
\]
Equation (17) corresponds to the equation with the Charlier polynomials as
eigenfunctions  (for the definition of the Meixner and Charlier polynomials
 see, e.g., \cite{nsu}). For a certain particular choice of the parameters,
one can reproduce the equations having Tschebyschov and Krawtchouk
polynomials as solutions. If $x$ is the continuous argument, we will arrive at
{\it continuous} analogues of above-mentioned polynomials. Up to our knowledge
those polynomials are not studied in literature.

Among the equations (13) there also exist quasi-exactly-solvable equations
possessing a finite number of polynomial eigenfunctions. All those equations
are classified via  the cubic polynomial element of the universal enveloping
$sl_2$-algebra taken in the representation (3), (8)
\begin{equation}
\label{e21}
\tilde T = A_+ (J^+_n + \delta J^0_nJ^0_n) + A_1 J^0_nJ^0_n
(J^-_n + {1 \over \delta}) + A_2 J^0_nJ^-_n + A_3 J^0_n + A_4 J^-_n + A_5
\end{equation}
(cf. (15)), where the $A$'s are free parameters.

In conclusion, it is worth emphasizing a quite surprising result:
 In general, three-point (quasi)-exactly-solvable {\it finite-difference}
operators of the type (13) emerging from (15), or (21) are isospectral
to (quasi)-exactly-solvable, third-order {\it differential} operators.

\end{document}